\documentstyle[twoside]{article}
\title{\bf Asymptotic expansions of Mathieu-Bessel series. I}
\author{\sc R. B. Paris\footnote{E-mail address:\ \ {\tt r.paris@abertay.ac.uk}}\\
\\
{\em Division of Computing and Mathematics,}\\
{\em Abertay University, Dundee DD1 1HG, UK}\\
}

\begin{document}
\newcommand{\bee}{\begin{equation}}
\newcommand{\ee}{\end{equation}}
\def\f#1#2{\mbox{${\textstyle \frac{#1}{#2}}$}}
\def\dfrac#1#2{\displaystyle{\frac{#1}{#2}}}
\newcommand{\fr}{\frac{1}{2}}
\newcommand{\fs}{\f{1}{2}}
\newcommand{\g}{\Gamma}
\newcommand{\om}{\omega}
\newcommand{\br}{\biggr}
\newcommand{\bl}{\biggl}
\newcommand{\ra}{\rightarrow}
\renewcommand{\topfraction}{0.9}
\renewcommand{\bottomfraction}{0.9}
\renewcommand{\textfraction}{0.05}
\newcommand{\mcol}{\multicolumn}
\newcommand{\gtwid}{\raisebox{-.8ex}{\mbox{$\stackrel{\textstyle >}{\sim}$}}}
\newcommand{\ltwid}{\raisebox{-.8ex}{\mbox{$\stackrel{\textstyle <}{\sim}$}}}
\date{}
\maketitle
\pagestyle{myheadings}
\markboth{\hfill {\it R.B. Paris} \hfill}
{\hfill {\it Asymptotics of a Mathieu-Bessel series  } \hfill}
\begin{abstract} 
We consider the asymptotic expansion of the Mathieu-Bessel series
\[S_\nu(a,b)=\sum_{n=1}^\infty \frac{n^\gamma J_\nu(nb/a)}{(n^2+a^2)^\mu}, \qquad (\mu, b>0,\ \gamma, \nu\in {\bf R})\]
as $a\to+\infty$ with the other parameters held fixed, where $J_\nu(x)$ is the Bessel function of the first kind of order $\nu$.
A special case arises when $\gamma+\nu$ is a positive even integer, where the expansion comprises finite algebraic terms together with an exponentially small expansion. Numerical examples are presented to illustrate the accuracy of the various expansions. The expansion of the alternating variant of $S_\nu(a,b)$ is considered. The series when the $J_\nu(x)$ function is replaced by the Bessel function $Y_\nu(x)$ is briefly mentioned.
\vspace{0.4cm}

\vspace{0.3cm}

\end{abstract}

\vspace{0.2cm}

\noindent $\,$\hrulefill $\,$

\vspace{0.2cm}

\begin{center}
{\bf 1. \  Introduction}
\end{center}
\setcounter{section}{1}
\setcounter{equation}{0}
\renewcommand{\theequation}{\arabic{section}.\arabic{equation}}
The functional series
\bee\label{e11}
\sum_{n=1}^\infty\frac{n}{(n^2+a^2)^\mu}
\ee
in the case $\mu=2$ 
was introduced by Mathieu in his 1890 book \cite{EM} dealing with the elasticity of solid bodies. Several integral representations for (\ref{e11}), and its alternating variant, have been obtained; see \cite{PST} and the references therein.

The asymptotic expansion of the more general Mathieu series
\bee\label{e12}
\sum_{n=1}^\infty\frac{n^\gamma}{(n^\beta+a^\beta)^\mu}\qquad (\mu>0,\ \beta>0,\ \beta\mu-\gamma>1)
\ee
was considered by Zvastavnyi \cite{VZ} for $a\ra+\infty$ and by Paris \cite{P13} for $|a|\ra\infty$ in the sector $|\arg\,a|<\pi/\beta$. In \cite{P13}, the additional factor $e_n:=\exp [-a^\beta b/(n^\beta+a^\beta)]$ (with $b>0$) was included in the summand which, although not affecting the rate of convergence of the series (since $e_n\ra1$ as $n\ra\infty$), can modify the large-$a$ growth, particularly with the alternating variant of (\ref{e12}). Both these authors adopted a Mellin transform approach. In \cite{P16} the above series was considered in the special case when $\beta$ and $\gamma$ assume even integer values. 
It was found that for these parameter values the asymptotic expansion for large complex $a$ in the sector $|\arg\,a|<\pi/\beta$  consists of a finite algebraic expansion together with an infinite sequence (when $\mu$ is not an integer) of increasingly subdominant exponentially small contributions. 

More recently, Gerhold and Tomovski \cite{GT} extended the asymptotic study of such Mathieu series by considering the power series
\[\sum_{n=1}^\infty \frac{nz^n}{(n^2+a^2)^\mu}\qquad (|z|\leq 1).\]
For $\mu>1$ and complex $z$ satisfying $|z|\leq 1$, with $z\neq 1$, they obtained the asymptotic expansion of this series as $a\to+\infty$ with coefficients involving the polylogarithm function $\mbox{Li}_\alpha(z)=\sum_{n\geq 1} n^{-\alpha} z^n$. From this result they were able to deduce the large-$a$ expansions of the trigonometric Mathieu series
\bee\label{e100}
\sum_{n=1}^\infty\frac{n \sin nx}{(n^2+a^2)^\mu},\qquad \sum_{n=1}^\infty\frac{n \cos nx}{(n^2+a^2)^\mu}.
\ee

In the present paper we extend (\ref{e100}) to the Mathieu-Bessel series
\bee\label{e101}
S_\nu(a,b)=\sum_{n=1}^\infty \frac{n^\gamma J_\nu(nb/a)}{(n^2+a^2)^\mu}\qquad (\mu>0,\ 2\mu-\gamma>\fs),
\ee
where $J_\nu(x)$ denotes the Bessel function of order $\nu$ and the condition $2\mu-\gamma>\fs$ secures absolute convergence of the series. The parameters will be supposed to satisfy $a>0$, $b>0$, $\mu>0$ with $\gamma, \nu\in{\bf R}$ and we shall be concerned with the limit $a\to\infty$. The alternating variant of this series is also considered. In the special case $\gamma+\nu=2m$, $m=0, 1, 2, \ldots\ $, the expansion of $S_\nu(a,b)$ is found to consist of an algebraic contribution together with an exponentially small expansion. This
mirrors the above-mentioned result found for the series in (\ref{e12}) when both $\beta$ and $\gamma$ are even integers. We present numerical results confirming the accuracy of the various expansions obtained and also briefly mention the series when the Bessel function is replaced by the Bessel function of the second kind $Y_\nu(x)$.

In the application of the Mellin transform method to the series in (\ref{e101}) and its alternating variant we shall require the following estimates for the gamma function and the Riemann zeta function. For real $\sigma$ and $t$, we have the estimates
\bee\label{e13}
\g(\sigma\pm it)=O(t^{\sigma-\frac{1}{2}}e^{-\frac{1}{2}\pi t}),\qquad |\zeta(\sigma\pm it)|=O(t^{\Omega(\sigma)} \log^\alpha t)\quad (t\ra+\infty),
\ee
where $\Omega(\sigma)=0$ ($\sigma>1$), $\fs-\fs\sigma$ ($0\leq\sigma\leq 1$), $\fs-\sigma$ ($\sigma<0$) and $\alpha=1$ ($0\leq\sigma\leq 1$), $\alpha=0$ otherwise \cite[p.~95]{ECT}. The zeta function $\zeta(s)$ has a simple pole of unit residue at $s=1$, $\zeta(0)=-\fs$ and has trivial zeros at $s=-2k$, $k=1, 2, \ldots\ $.
Finally, we have the well-known functional relation satisfied by $\zeta(s)$ given by \cite[p.~603]{DLMF}
\bee\label{e15}
\zeta(s)=2^s \pi^{s-1} \zeta(1-s) \g(1-s) \sin \fs\pi s.
\ee
\vspace{0.6cm}

\begin{center}
{\bf 2. \ The asymptotic expansion of $S_\nu(a,b)$ for $a\to+\infty$}
\end{center}
\setcounter{section}{2}
\setcounter{equation}{0}
\renewcommand{\theequation}{\arabic{section}.\arabic{equation}}
Let $a>0$, $b>0$, $\mu>0$ and $\gamma$, $\nu$ be real parameters.
We consider the asymptotic expansion of the Mathieu-Bessel series 
\bee\label{e21}
S_\nu(a,b)=\sum_{n=1}^\infty \frac{n^\gamma J_\nu(nb/a)}{(n^2+a^2)^\mu}\qquad (2\mu-\gamma>\fs)
\ee
as $a\to+\infty$. The above sum may be  written as
\[S_\nu(a,b)=a^{\gamma-2\mu}\sum_{n=1}^\infty h(n/a),\qquad
h(x):=\frac{x^\gamma J_\nu(bx)}{(1+x^2)^\mu}.\]

We employ a Mellin transform approach as discussed in, for example, \cite[Section 4.1.1]{PK}. The Mellin transform of $h(x)$ is $H(s)=\int_0^\infty x^{s-1}h(x)\,dx$, and with
\[\lambda\equiv\lambda(s):=\frac{1}{2}(s+\gamma+\nu),\qquad\chi:=\frac{1}{4}b^2, \qquad B:=\frac{(b/2)^\nu}{2\g(\mu)},\]
we have the evaluation \cite[p.~434]{WBF}
\[H(s)=\frac{\pi B}{\sin \pi(\mu\!-\!\lambda)}\bl\{\g(\lambda) {}_1{\bf F}_2\bl(\!\!\begin{array}{c}\lambda\\1\!+\!\lambda\!-\!\mu, 1\!+\!\nu\end{array}\bl|\,\chi\br)\hspace{5cm}\]
\bee\label{e22}
\hspace{4cm}-\chi^{\mu-\lambda} \g(\mu){}_1{\bf F}_2\bl(\!\!\begin{array}{c}\mu\\1\!-\!\lambda\!+\!\mu, 1\!-\!\lambda\!+\!\mu\!+\!\nu\end{array}\bl|\,\chi\br)\br\}
\ee
valid in the strip
\[-\gamma-\nu<\Re (s)<\delta,\qquad \delta:=2\mu-\gamma+\fs.\]
We note that $\delta>1$ on account of the convergence condition in (\ref{e21}).
The hypergeometric functions appearing in (\ref{e22}) are defined by
\[{}_1{\bf F}_2\bl(\!\!\begin{array}{c}\alpha\\ \beta, \gamma\end{array}\bl|z\br)=\frac{1}{\g(\beta) \g(\gamma)}\,{}_1 F_2\bl(\!\!\begin{array}{c}\alpha\\ \beta, \gamma\end{array}\bl|z\br)=\frac{1}{\g(\beta) \g(\gamma)}\,\sum_{n=0}^\infty \frac{(\alpha)_n}{(\beta)_n (\gamma)_n}\,\frac{z^n}{n!},\]
where $(\alpha)_n=\g(\alpha+n)/\g(\alpha)$ is the Pochhammer symbol.
Use of the Mellin inversion theorem (see, for example, \cite[p.~118]{PK}) then leads to
\bee\label{e23}
S_\nu(a,b)=\frac{a^{\gamma-2\mu}}{2\pi i} \sum_{n=1}^\infty\int_{c-\infty i}^{c+\infty i}H(s) (n/a)^{-s}ds=
\frac{a^{\gamma-2\mu}}{2\pi i} \int_{c-\infty i}^{c+\infty i}H(s)\zeta(s)a^s\,ds,
\ee
where $\max\{-\gamma-\nu,1\}<c<\delta$ and $\zeta(s)$ denotes the Riemann zeta function. The inversion of the order of summation and integration is justified by absolute convergence provided $c$ satisfies the stated condition.

We now need to examine the pole structure of $H(s)$. The function ${}_1{\bf F}_2(z)$ is an entire function of its parameters (and hence of the variable $s$) \cite[(16.2.5)]{DLMF} and consequently has no poles. The poles of $H(s)$ are situated at $\lambda=-k$ and $\mu-\lambda=\pm k$, $k=0, 1, 2, \ldots\ $. However, it is shown in Appendix A that $H(s)$ is in fact regular at the points $\mu-\lambda=\pm k$ due to a cancellation of terms. Thus, the only singularities of $H(s)$ are those corresponding to $\lambda=-k$; that is, at $s=-\gamma-\nu-2k$, $k=0, 1, 2 \ldots\ $. In addition, there is a pole of the integrand in (\ref{e23}) at $s=1$ resulting from $\zeta(s)$.

We consider the integral in (\ref{e23}) taken round the rectangular contour with vertices at $c\pm iT$ and $-c'\pm iT$,
where $c'>0$ and $T>0$. Then the leading behaviour of the first and second hypergeometric functions in (\ref{e22}) can be shown to
be $(\fs b)^{-\nu}I_\nu(b)/\g(1+\lambda-\mu)$, where $I_\nu(b)$ is the modified Bessel function, and \[\frac{O(1)}{\g(1\!-\!\lambda\!+\!\mu)\g(1\!-\!\lambda\!+\!\mu\!+\!\nu)}=\frac{O(1)\g(\lambda-\mu)}{\g(1-\lambda+\mu+\nu)}
\,\sin \pi(\lambda-\mu)\]
respectively for large $\lambda$. We therefore find on the upper and lower sides of the rectangle $s=\sigma\pm iT$, $-c'\leq\sigma\leq c$, with $\lambda=\fs(\sigma+\gamma+\nu)\pm\fs iT$, that
\[|H(s)|=O\bl(\frac{\g(\lambda)}{\g(1+\lambda-\mu) \sin \pi(\mu-\lambda)}\br)+\frac{O(1)\g(\lambda-\mu)}{\g(1-\lambda+\mu+\nu)}\]
\[\hspace{0.2cm}=O(T^{\mu-1}e^{-\pi T/2})+O(T^{\sigma-\delta-\frac{1}{2}})=O(T^{\sigma-\delta-\frac{1}{2}})\]
as $T\to\infty$ by the first formula in (\ref{e13}).
Thus, the modulus of the integrand has on these paths has the order estimate 
\[|H(s)a^s \zeta(s)|=O(T^{\omega(\sigma)} \log\,T),\qquad\omega(\sigma):=\Omega(\sigma)+\sigma-\delta-\fs,\]
where, from (\ref{e13}), $\omega(\sigma)=\sigma-\delta-\fs$ ($\sigma>1$), $\fs\sigma-\delta$ ($0\leq\sigma\leq1$), $-\delta$ ($\sigma<0$). Since $\sigma<\delta$ and $\delta>1$, the contribution from the upper and lower sides of the rectangle therefore vanishes as $T\to\infty$.

Displacement of the integration path to the left over the simple poles at $s=1$ and $s=-\gamma-\nu-2k$ (when $\gamma+\nu$ is not a negative odd integer) then yields the following result.
\newtheorem{theorem}{Theorem}
\begin{theorem}$\!\!\!.$ \ Let  $\chi=b^2/4$ and $\gamma+\nu\neq -1, -3, \ldots\ $. Then we have the asymptotic expansion
\bee\label{e24}
S_\nu(a,b)-a^{\gamma-2\mu+1} H(1)\sim \frac{a^{-\nu-2\mu}(b/2)^\nu}{\g(1+\nu)} \sum_{k=0}^\infty\frac{(-)^k}{k!}\,(\mu)_k\zeta(-\omega_k) F_k\,a^{-2k}
\ee
as $a\to+\infty$, where $\omega_k:=\gamma+\nu+2k$ and
\bee\label{e24a}
F_k:={}_1F_2\bl(\!\!\begin{array}{c}-k\\ 1\!-\!\mu\!-\!k, 1\!+\!\nu\end{array}\bl| \chi\br).
\ee
\end{theorem}

The functional relation for $\zeta(s)$ in (\ref{e15}) can be employed to obtain an alternative version of the asymptotic series on the right-hand side of (\ref{e24}) in the form
\[\hspace{2cm}-\frac{a^{\nu-2\mu}(b/2)^\nu \sin \fs\pi(\gamma+\nu)}{2^{\gamma+\nu} \pi^{1+\gamma+\nu} \g(1+\nu)} \sum_{k=0}^\infty \frac{(\mu)_k}{k!} \g(1+\omega_k)\zeta(1+\omega_k)F_k (2\pi a)^{-2k}.\] 
The $F_k$ are polynomials in $\chi$ of degree $k$ and
\[F_0=1,\quad F_1=1+\frac{\chi}{\mu(1+\nu)},\quad F_2=1+\frac{2\chi}{(1+\mu)(1+\nu)}+\frac{\chi^2}{\mu(1+\mu)(1+\nu)_2}.\]
In the special case $\mu=1$, the $F_k$ are independent of $k$ and given by
\bee\label{e24b}
F_k={}_0F_1\bl(\!\!\begin{array}{c} -\!\!-\\1+\nu\end{array}\bl|\chi\br)=\g(1+\nu) (\fs b)^{-\nu} I_\nu(b)\equiv {\cal I}_\nu(b).
\ee

The case $\gamma+\nu=2m$, $m=0, 1, 2, \ldots$ requires separate attention, which is discussed in the next section.
When $\gamma+\nu=-1, -3, \ldots$ the pole at $s=1$ is double and the residue has to be evaluated accordingly. We present only the case $\gamma+\nu=-1$; other cases can be dealt with in a similar manner. In Appendix B, it is established that the residue at $s=1$ when $\gamma+\nu=-1$, $\mu\neq 1, 2, \ldots$ is given by
\bee\label{e25}
\mbox{Res}_{s=1}=\frac{2aB\g(\mu)}{\g(1+\nu)}\bl\{\log\,a+\frac{1}{2}(A+{\hat \gamma}-\psi(\mu))\br\},
\ee
where ${\hat \gamma}=0.57721\ldots$ is the Euler-Mascheroni constant, $\psi(x)$ is the logarithmic derivative of $\g(x)$ and
\[A=\g(1+\nu)\bl\{\chi \g(1-\mu)\,{}_2{\bf F}_3\bl(\!\!\begin{array}{c}1,1\\2,2\!-\!\mu, 2\!+\!\nu\end{array}\bl|\chi\br)-\frac{\pi\chi^\mu}{\sin \pi\mu}\,{}_1{\bf F}_2\bl(\!\!\begin{array}{c}\mu\\1\!+\!\mu, 1\!+\!\mu\!+\!\nu\end{array}\bl|\chi\br)\br\}.\]
In the special case $\mu=1$ (and $\gamma+\nu=-1$) we have by a limiting process (see (\ref{b2}))
\bee\label{e26}
\mbox{Res}_{s=1}=\frac{aB}{\g(1+\nu)}\bl\{2\log\,a+2{\hat\gamma}+\kappa(1-{\cal I}_\nu(b))-\frac{\chi F_*}{1+\nu}\br\}
\ee
where $\kappa:=1-{\hat\gamma}+\psi(2+\nu)-\log\,\chi$,
${\cal I}_\nu(b)$ is defined in (\ref{e24b}) and $F_*$ is defined in (\ref{b3}).
%
This produces the result:
\begin{theorem}$\!\!\!.$ \ Let $\chi=b^2/4$ and $B=(b/2)^\nu/(2\g(\mu))$. Then, when $\gamma+\nu=-1$, we have the asymptotic expansion
\bee\label{e27}
S_\nu(a,b)-\mbox{Res}_{s=1}\sim\sum_{k=1}^\infty\frac{(-)^k}{k!}\,(\mu)_k\zeta(1\!-\!2k) F_k\,a^{-2k}
\ee
as $a\to+\infty$, where $\mbox{Res}_{s=1}$ is given by (\ref{e25}) when $\mu\neq 1, 2, \ldots$ and by (\ref{e26}) in the special case $\mu=1$. The quantities $F_k$ are specified in (\ref{e24a}).
\end{theorem}
\vspace{0.6cm}
 
\begin{center}
{\bf 3. \ The exponentially small case when $\gamma+\nu$ is a non-negative even integer}
\end{center}
\setcounter{section}{3}
\setcounter{equation}{0}
\renewcommand{\theequation}{\arabic{section}.\arabic{equation}}
In this section we let $\gamma+\nu=2m$, $m=0, 1, 2, \ldots\ $. Then the asymptotic series on the right-hand side of (\ref{e24}) vanishes
when $m\geq 1$, and consists of a single term when $m=0$, on account of the trivial zeros of $\zeta(s)$. In this case the expansion of $S_\nu(a,b)$ will involve an exponentially small contribution as $a\to+\infty$; a similar result was obtained in \cite{P16} for the Mathieu series (\ref{e12}) when both $\beta$ and $\gamma$ are even integers.

To deal with this case, we first observe that the only poles of $H(s)$ are situated at $s=-2m-2k$, $k=0, 1, 2, \ldots$
which, apart from the case $m=0$, are all cancelled by the trivial zeros of $\zeta(s)$. When $m=0$, there is a pole at $s=0$ resulting from the first part of the expression in (\ref{e22}). Then, displacement of the integration path in (\ref{e23}) to the left into $\Re (s)<0$ over the poles at $s=1$ and $s=0$ (when $m=0$) yields
\[S_\nu(a,b)-a^{\gamma-2\mu+1}H(1)+\frac{a^{-\nu-2\mu}(b/2)^\nu}{2\g(1+\nu)}\,\delta_{0m}=\frac{a^{\gamma-2\mu}}{2\pi i} \int_{-c-\infty i}^{-c+\infty i} H(s) \zeta(s)a^sds\qquad(c>0),\]
where $\delta_{0m}$ denotes the Kronecker delta symbol. We now
make the change of variable $s\to -s$ and use the functional relation for $\zeta(s)$ in (\ref{e15}) to obtain
\[S_\nu(a,b)-a^{\gamma-2\mu+1}H(1)+\frac{a^{-\nu-2\mu}(b/2)^\nu}{2\g(1+\nu)}\,\delta_{0m}\hspace{5cm}\]
\bee\label{e30}
\hspace{4cm}=-\frac{a^{\gamma-2\mu}}{2\pi i} \int_L H(-s) \zeta(1+s) \g(1+s) \frac{\sin \fs\pi s}{\pi}\,(2\pi a)^{-s}ds,
\ee
where the contour $L$ denotes a path parallel to the imaginary axis that can be displaced as far to the right as we please since there are no poles of the integrand in $\Re (s)>0$. On the path $L$ we set $s=\sigma+it$, $t\in (-\infty, \infty)$ and choose $\sigma=2N+1-2\mu$, where $N$ denotes a (large) positive integer.

The function $H(-s)$ may now be split into two parts: $H_1(-s)$ and $H_2(-s)$, where $H_1(-s)$ corresponds to the part containing the first ${}_1{\bf F}_2(\chi)$ function in (\ref{e22}) and $H_2(-s)$ to the second such function. Then, from (\ref{e22}) with $\lambda=m-\fs s$, we obtain
\bee\label{e31}
H_1(-s) \g(1+s)\,\frac{\sin \fs\pi s}{\pi}=\frac{(-)^{m-1}B}{\g(1+\nu)}\,G_1(s)\,{}_1F_2\bl(\!\!\begin{array}{c}\!\!m\!-\!\fs s\\1\!-\!\mu\!+\!m\!-\!\fs s, 1\!+\!\nu\end{array}\bl|\,\chi\br),
\ee
where
\[G_1(s):=\frac{\g(1+s) \g(\mu-m+\fs s)}{\g(1-m+\fs s)}.\]
Since $|s|$ is everywhere large on the contour $L$, the ratio of gamma functions appearing in (\ref{e31}) may be expanded as an inverse factorial expansion in the form \cite[p.~53]{PK} 
\bee\label{e32}
G_1(s)\sim 2^{1-\mu} \g(s+\mu)\bl\{1+\frac{C_1}{s\!+\!\mu\!-\!1}+\frac{C_2}{(s\!+\!\mu\!-\!1)(s\!+\!\mu\!-\!2)}+\cdots\br\},
\ee
where, with some effort, it may be shown that
\[C_1=\frac{1}{2}(1-\mu)(4m-\mu),\qquad C_2=\frac{1}{8}(1-\mu)(2-\mu)\{16m^2+(1+\mu)(\mu-8m)\}, \ldots\ .\]

A similar type of expansion for the hypergeometric function appearing in (\ref{e31}) is considered in Appendix C. It is shown that
\bee\label{e33}
{}_1F_2\bl(\!\!\begin{array}{c}m\!-\!\fs s\\1\!-\!\mu\!+\!m\!-\!\fs s, 1\!+\!\nu\end{array}\bl|\,\chi\br)=\g(1+\nu) (\fs b)^{-\nu} I_\nu(b) \bl\{1+\frac{C_1'}{s\!+\!\mu\!-\!1}+\frac{C_2'}{(s\!+\!\mu\!-\!1)(s\!+\!\mu\!-\!2)}+\cdots\br\}
\ee
as $s\to\infty$ in $\Re (s)>0$, where
\[C_1'=(1-\mu)b\,\frac{I_{\nu+1}(b)}{I_\nu(b)},\quad C_2'=(1-\mu)b\bl\{(2m\!+\!1\!-\!\mu)\,\frac{I_{\nu+1}(b)}{I_\nu(b)}+\frac{(2-\mu)b}{2}\,\frac{I_{\nu+2}(b)}{I_\nu(b)}\br\}.\]
It then follows from (\ref{e32}) and (\ref{e33}) that
\bee\label{e34}
G_1(s){}_1F_2\bl(\!\!\begin{array}{c}\!\!m-\fs s\\1\!-\!\mu\!+\!m\!-\!\fs s, 1\!+\!\nu\end{array}\bl|\,\chi\br)
=2^{1-\mu} \g(1+\nu) (\fs b)^{-\nu} I_\nu(b) \sum_{j\geq0} D_j \g(s+\mu-j),
\ee
where
\[D_0=1,\quad D_1=C_1+C_1'=(1-\mu)\bl\{\frac{b I_{\nu+1}(b)}{I_\nu(b)}+2m-\frac{1}{2}\mu\br\},\]
\[D_2=C_2+C_2'+C_1C_1'\]
\[=(1-\mu)b\bl\{[2m+(1-\mu)(2m+1)+\frac{1}{2}\mu(1-\mu)]\,\frac{I_{\nu+1}(b)}{I_\nu(b)}+\frac{(2-\mu)b}{2}\,\frac{I_{\nu+2}(b)}{I_\nu(b)}\br\}\]
\bee\label{e3d}
+\frac{1}{8}(1-\mu)(2-\mu)\{16m^2+(\mu+1)(\mu-8m)\}, \ldots\ .
\ee
Insertion of the expansion (\ref{e34}) into the integral on the right-hand side of (\ref{e30}), with $H(-s)$ replaced by $H_1(-s)$,
then produces (formally)
\[S_\nu^{(1)}(a,b):=-\frac{a^{\gamma-2\mu}}{2\pi i} \int_L H_1(-s) \zeta(1+s) \g(1+s) \frac{\sin \fs\pi s}{\pi}\,(2\pi a)^{-s}ds\]
\[=\frac{(-)^m a^{\gamma-2\mu}}{2^{\mu}\g(\mu)}\,I_\nu(b)\,\sum_{j\geq0} \frac{D_j}{2\pi i} \int_L \zeta(1+s)\g(s+\mu-j)\,(2\pi a)^{-s}ds\]
\[=\frac{(-)^m a^{\gamma-2\mu}}{2^{\mu}\g(\mu)}\,I_\nu(b)\,\sum_{n\geq 1}\frac{1}{n}\sum_{j\geq0} \frac{D_j}{2\pi i} \int_L \g(s+\mu-j)\,(2\pi na)^{-s}ds\]
upon expansion of $\zeta(1+s)$ into its power series representation (which is permissible since $\Re (s)>1$ on $L$). 

The above integrals can be evaluated by means of the Cahen-Mellin integral \cite[p.~90]{PK}
\[\frac{1}{2\pi i}\int_{c-\infty i}^{c+\infty i} \g(s+\alpha) z^{-s}ds=z^{\alpha} e^{-z} \qquad(c>0,\ |\arg\,z|<\fs\pi)\]
to yield the expansion
\[S_\nu^{(1)}(a,b)\sim (-)^m \frac{a^{\gamma-\mu} \pi^\mu}{\g(\mu)}\,I_\nu(b) \sum_{n\geq 1}\frac{e^{-2\pi na}}{n^{1-\mu}} \sum_{j\geq0}\frac{D_j}{(2\pi na)^j}\]
as $a\to+\infty$. Thus we find the dominant large-$a$ contribution when $\gamma+\nu=2m$ given by
\bee\label{e35}
S_\nu^{(1)}(a,b)\sim(-)^m \frac{a^{\gamma-\mu} \pi^\mu}{\g(\mu)}\,I_\nu(b)\bl\{e^{-2\pi a} \sum_{j\geq0}\frac{D_j}{(2\pi a)^j}+O(e^{-4\pi a})\br\}.
\ee

The contribution to the integral in (\ref{e30}) resulting from $H_2(-s)$ is
\[S_\nu^{(2)}(a,b):=-\frac{a^{\gamma-2\mu}}{2\pi i} \int_L H_2(-s) \zeta(1+s) \g(1+s) \frac{\sin \fs\pi s}{\pi}\,(2\pi a)^{-s}ds,\]
where, from (\ref{e22}),
\[H_2(-s)=\frac{(-)^{m-1}\pi B \g(\mu)\chi^{\mu-m+s/2}\csc \pi(\mu+\fs s)}{\g(1\!+\!\mu\!-\!m\!+\!\fs s)\g(1\!+\!\mu\!+\!\nu\!-\!m\!+\!\fs s)}\hspace{3cm}\]
\[\hspace{5cm}\times {}_1F_2\bl(\!\!\begin{array}{c}\mu\\1\!+\!\mu\!-\!m\!+\!\fs s, 1\!+\!\mu\!+\!\nu\!-\!m\!+\!\fs s\end{array}\bl|\,\chi\br).\]
It is readily seen that the above hypergeometric function satisfies ${}_1F_2(\chi)\to 1$ as $s\to\infty$ in $\Re (s)>0$. Then the quantity
\[H_2(-s) \g(1+s)\,\frac{\sin \fs\pi s}{\pi}\sim (-)^{m-1} B\g(\mu) \chi^{\mu-m+s/2} G_2(s)\,\frac{\sin \fs\pi s}{\sin \pi(\mu+\fs s)},\]
where 
\[G_2(s):=\frac{\g(1+s)}{\g(1\!+\!\mu\!-\!m\!+\!\fs s)\g(1\!+\!\mu\!+\!\nu\!-\!m\!+\!\fs s)}\]
\[\hspace{2.6cm}\sim \pi^{-1/2} 2^{s-\alpha+1/2} s^{\alpha-1},\qquad \alpha:=2m-2\mu-\nu+\fs\]
as $s\to\infty$. On $L$ we have (with $s=\sigma+it$, $\sigma=2N+1-2\mu$)
\[\bl|\frac{\sin \fs\pi s}{\sin \pi(\mu+\fs s)}\br|=\bl(\frac{\sin^2\fs\pi\sigma+\sinh^2\fs\pi t}{1+\sinh^2 \fs\pi t}\br)^{\!1/2}\leq 1.\]

Hence, using the result that $|\zeta(x+iy)|<\zeta(x)$ for real $x>1$ and $y$, we find that
\[|S_\nu^{(2)}(a,b)|<\frac{a^{\gamma-2\mu}B\g(\mu)\zeta(N+1)}{2^{\alpha+1/2} \pi^{3/2}} \int_LX^{-s} s^{\alpha-1}ds,\qquad X:=\frac{2\pi a}{b}.\]
Now
\[\int_LX^{-s} s^{\alpha-1}ds=\int_L s^{\alpha-1} e^{-s \log\,X}\,ds=-\bl[\frac{\g(\alpha, s\log\,X)}{(\log X)^\alpha}\br]_{\!L},\]
where $\g(a,z)$ denotes the incomplete gamma function. Since \cite[(8.11.2)]{DLMF} 
\[\g(\alpha,s\log X)\sim (s\log X)^{\alpha-1} X^{-s}\qquad (|s|\log X\to\infty,\ X>1,\  \Re (s)>0),\]
we see that this last integral vanishes as $N\to\infty$. Hence it follows that $S_\nu^{(2)}(a,b)\to 0$ as $N\to\infty$. 

Then we obtain the following expansion:
\begin{theorem}$\!\!\!.$\ When $\gamma+\nu=2m$, $m=0, 1, 2, \ldots$ we have the exponentially small asymptotic expansion
\[S_\nu(a,b)-a^{\gamma-2\mu+1} H(1)+\frac{a^{-\nu-2\mu}(b/2)^\nu}{2\g(1+\nu)}\,\delta_{0m}\hspace{5cm}\]
\bee\label{e38}
\hspace{4cm}\sim(-)^m \frac{a^{\gamma-\mu} \pi^\mu}{\g(\mu)}\,I_\nu(b)\bl\{e^{-2\pi a} \sum_{j\geq0}\frac{D_j}{(2\pi a)^j}+O(e^{-4\pi a})\br\}
\ee
as $a\to+\infty$, where the first three coefficients $D_j$ are given in (\ref{e3d}) and $\delta_{0m}$ is the Kronecker delta symbol. The function $H(s)$ is defined in (\ref{e22}).
\end{theorem}
\vspace{0.6cm}

\begin{center}
{\bf 4. \ Numerical results and concluding remarks}
\end{center}
\setcounter{section}{4}
\setcounter{equation}{0}
\renewcommand{\theequation}{\arabic{section}.\arabic{equation}}
We present some numerical results to illustrate the accuracy of the various expansions obtained. The values shown give the absolute relative error in the high-precision computation of the Mathieu-Bessel series from (\ref{e101}) using the asymptotic expansions in Theorems 1 and 2.  In Table 1 we show values\footnote{In the tables we write the values as $x(y)$ instead of $x\times 10^y$.} of
\bee\label{e41}
{\hat S}:=-S_\nu(a,b)+a^{\gamma-2\mu+1}H(1)
\ee
for different values of $a$ and the absolute relative error resulting from optimal truncation of (\ref{e24}) (that is, truncation just before the least term in the expansion). In Table 2 we show the values of ${\tilde S}:=S_\nu(a,b)-\mbox{Res}_{s=1}$, where $\mbox{Res}_{s=1}$ is defined in (\ref{e25}) (and (\ref{e26}) when $\mu=1$), for different values of $a$ in the double-pole case $\gamma+\nu=-1$. The associated absolute relative error resulting from optimal truncation of the asymptotic series in (\ref{e27}) is indicated.
\begin{table}[t]
\caption{\footnotesize{The absolute relative error in the computation of ${\hat S}$ in (\ref{e41}) from (\ref{e24}) for different $a$ and $b$  when $\mu=3$, $\nu=1/3$ and $\gamma=1/2$. The optimal truncation index $k_o$ employed in the asymptotic sum (\ref{e24}) is indicated.}}
\begin{center}
\begin{tabular}{|rc|cc|cc|cc|}
\hline
&&&&&&&\\[-0.25cm]
\mcol{2}{|c|}{} & \mcol{2}{c|}{$b=\fs$} &\mcol{2}{c|}{$b=1$}& \mcol{2}{c|}{$b=2$} \\
\mcol{1}{|c}{$a$} & \mcol{1}{c|}{$k_o$} & \mcol{1}{c}{${\hat S}$} & \mcol{1}{c|}{Error}& \mcol{1}{c}{${\hat S}$} & \mcol{1}{c|}{Error} & \mcol{1}{c}{${\hat S}$} & \mcol{1}{c|}{Error}\\
[.1cm]\hline
&&&&&&&\\[-0.25cm]
2 & 4 &  $1.08383(-03)$ & $\!\!\!\!1.867(-03)$ & $1.37088(-03)$ & $\!\!\!\!2.107(-02)$ & $1.75474(-03)$ & $\!\!\!\!3.236(-03)$ \\
4 & 10 & $1.26425(-05)$ & $\!\!\!\!6.507(-08)$ & $1.59404(-05)$ & $\!\!\!\!7.411(-08)$ & $2.01438(-05)$ & $\!\!\!\!1.176(-07)$ \\
6 & 16 & $9.60936(-07)$ & $\!\!\!\!2.034(-12)$ & $1.21109(-06)$ & $\!\!\!\!2.800(-12)$ & $1.52781(-06)$ & $\!\!\!\!1.361(-12)$ \\
8 & 22 & $1.54920(-07)$ & $\!\!\!\!5.349(-13)$ & $1.95222(-07)$ & $\!\!\!\!4.492(-14)$ & $2.46137(-07)$ & $\!\!\!\!5.286(-14)$ \\
[.1cm]\hline
\end{tabular}
\end{center}
\end{table}
\begin{table}[t]
\caption{\footnotesize{The absolute relative error in the computation of ${\tilde S}\equiv S_\nu(a,b)-\mbox{Res}_{s=1}$ for different $a$ when $\gamma+\nu=-1$ and $b=1$. The optimal truncation index $k_o$ employed in the asymptotic sum (\ref{e27}) is indicated.}}
\begin{center}
\begin{tabular}{|rc|cc|cc|cc|}
\hline
&&&&&&&\\[-0.25cm]
\mcol{2}{|c|}{} & \mcol{2}{c|}{$\gamma=-1,\ \nu=0;\ \mu=\f{5}{2}$} &\mcol{2}{c|}{$\ \gamma=-\f{3}{4},\ \nu=-\f{1}{4};\ \mu=\f{5}{2}$}& \mcol{2}{c|}{$\gamma=-\f{9}{4},\ \nu=\f{5}{4};\ \mu=1$} \\
\mcol{1}{|c}{$a$} & \mcol{1}{c|}{$k_o$} & \mcol{1}{c}{${\tilde S}$} & \mcol{1}{c|}{Error}& \mcol{1}{c}{${\tilde S}$} & \mcol{1}{c|}{Error} & \mcol{1}{c}{${\tilde S}$} & \mcol{1}{c|}{Error}\\
[.1cm]\hline
&&&&&&&\\[-0.25cm]
2 & 5 &  $1.89563(-03)$ & $\!\!\!\!3.209(-04)$ & $2.25578(-03)$ & $\!\!\!\!3.212(-04)$ & $9.29534(-04)$ & $\!\!\!\!3.088(-05)$ \\
4 & 11 & $1.41547(-05)$ & $\!\!\!\!6.178(-08)$ & $2.00180(-05)$ & $\!\!\!\!1.841(-08)$ & $2.38806(-05)$ & $\!\!\!\!3.446(-10)$ \\
6 & 17 & $8.22875(-07)$ & $\!\!\!\!8.322(-13)$ & $1.28777(-06)$ & $\!\!\!\!3.447(-13)$ & $2.83132(-06)$ & $\!\!\!\!3.522(-13)$ \\
8 & 24 & $1.09590(-07)$ & $\!\!\!\!9.066(-15)$ & $1.84289(-07)$ & $\!\!\!\!5.132(-14)$ & $6.24484(-07)$ & $\!\!\!\!3.276(-14)$ \\
[.1cm]\hline
\end{tabular}
\end{center}
\end{table}

The exponentially small case is illustrated in Table 3 where we show values of the absolute relative error in the computation of
\bee\label{e42}
{\cal S}:=S_\nu(a,b)-a^{\gamma-2\mu+1}H(1)+\frac{a^{-\nu-2\mu}(b/2)^\nu}{2\g(1+\nu)}\,\delta_{0m}
\ee
in the case $\gamma+\nu=0$ ($m=0$) by means of the expansion (\ref{e38}) as a function of truncation index$j$.
\begin{table}[h]
\caption{\footnotesize{The absolute relative error in the computation of ${\cal S}$ in (\ref{e42}) using the expansion in (\ref{e38}) as a function of truncation index $j\leq 2$ when $\gamma+\nu=0$ and $a=8$, $\mu=4$.}}
\begin{center}
\begin{tabular}{|r|c|c|c|}
\hline
&&&\\[-0.25cm]
\mcol{1}{|c|}{} & \mcol{1}{c|}{$\gamma=0,\ \nu=0$} &\mcol{1}{c|}{$\gamma=0,\ \nu=0$}& \mcol{1}{c|}{$\gamma=-\f{1}{3}, \nu=\f{1}{3}$} \\
\mcol{1}{|c|}{$j$} & \mcol{1}{c|}{$b=1$} & \mcol{1}{c|}{$b=3$} & \mcol{1}{c|}{$b=1$}\\
[.1cm]\hline
&&&\\[-0.25cm]
0 &  $9.729(-02)$ & $2.346(-02)$ & $1.035(-01)$ \\
1 &  $4.176(-03)$ & $2.257(-03)$ & $4.441(-03)$ \\
2 &  $8.129(-05)$ & $1.959(-06)$ & $8.621(-05)$ \\
[.1cm]\hline
\end{tabular}
\end{center}
\end{table}

When $\nu=\fs$, the Bessel function $J_\nu(z)$ reduces to $(2/\pi z)^{1/2} \sin z$, so that
\[S_{1/2}(a,b)=\bl(\frac{2a}{\pi b}\br)^{\!1/2} \sum_{n=1}^\infty \frac{n^{\gamma-1/2}}{(n^2+a^2)^\mu}\,\sin nx,\qquad x:=\frac{b}{a}.\]
An analogous result holds when $\nu=-\fs$, with the trigonometric function then replaced by $\cos nx$. Although these cases bear a superficial similarity to the series considered by Gerhold and Tomovski stated in (\ref{e100}) (see also \cite[Section 4]{GT}), the asymptotic expansions are different. These authors considered the limit $a\to\infty$ with $x$ fixed, whereas our expansion corresponds to $a\to\infty$, $x\to 0$.

We remark that the asymptotic expansion of the alternating version of (\ref{e101}) can be deduced by making use of the identity
\bee\label{e43}
\sum_{n=1}^\infty \frac{(-)^{n-1} n^\gamma }{(n^2+a^2)^\mu}\,J_\nu(nb/a)=S_\nu(a,b)-2^{\gamma-2\mu+1} S_\nu(\fs a,b).
\ee
From (\ref{e24}), this then produces the asymptotic expansion (when $\gamma+\nu\neq -1, -3, \ldots$)
\[\sum_{n=1}^\infty \frac{(-)^{n-1} n^\gamma }{(n^2+a^2)^\mu}\,J_\nu(nb/a)\sim
\frac{a^{\nu-2\mu}(b/2)^\nu}{\g(1+\nu)} \sum_{k=0}^\infty \frac{(-)^k(\mu)_k}{k!}\zeta(-\omega_k)F_k \{1-2^{1+\omega_k}\}a^{-2k}\]
as $a\to\infty$, where we recall that $\omega_k=\gamma+\nu+2k$.

Finally, the series involving the $Y$-Bessel function can, in the case of non-integer $\nu$, be obtained from the standard definition of $Y_\nu(z)$ in terms of $J_{\pm\nu}(z)$ given by
$Y_\nu(z)=\cot \pi\nu\,J_\nu(z)-\csc \pi\nu\,J_{-\nu}(z)$, to yield
\bee\label{e44}
\sum_{n=1}^\infty \frac{n^\gamma Y_\nu(nb/a)}{(n^2+a^2)^\mu}=\cot \pi\nu\, S_\nu(a,b)-\csc \pi\nu\, S_{-\nu}(a,b).
\ee
The case of integer values of $\nu$ would require a limiting procedure, which we do not consider further here.
\vspace{0.6cm}

\begin{center}
{\bf Appendix A: Discussion of the pole structure of $H(s)$}
\end{center}
\setcounter{section}{1}
\setcounter{equation}{0}
\renewcommand{\theequation}{\Alph{section}.\arabic{equation}}\vspace{0.6cm}
We examine the pole structure of the function $H(s)$ defined in (\ref{e22}), which has apparent singularities when $\mu-\lambda=\pm k$, $k=0, 1, 2, \ldots\ $, where $\lambda=\fs(\gamma+\nu+s)$. To demonstrate that $H(s)$ is regular at these points, it will be sufficient to show that 
\[Q(s):=\g(\lambda) {}_1{\bf F}_2\bl(\!\!\begin{array}{c}\lambda\\1\!+\!\lambda\!-\!\mu, 1\!+\!\nu\end{array}\bl|\,\chi\br)
- \g(\mu)\chi^{\mu-\lambda}{}_1{\bf F}_2\bl(\!\!\begin{array}{c}\mu\\1\!-\!\lambda\!+\!\mu, 1\!-\!\lambda\!+\!\mu\!+\!\nu\end{array}\bl|\,\chi\br)\]
has simple zeros when $\mu-\lambda=\pm k$.

First, consider the case $\mu-\lambda=k$, with $k=0, 1, 2, \ldots\ $. Then
\bee\label{a1}
Q(s)=\g(\lambda) {}_1{\bf F}_2\bl(\!\!\begin{array}{c}\mu\!-\!k\\1\!-\!k, 1\!+\!\nu\end{array}\bl|\,\chi\br)
- \g(\mu)\chi^{k}{}_1{\bf F}_2\bl(\!\!\begin{array}{c}\mu\\k\!+\!1, 1\!+\!\nu\!+\!k\end{array}\bl|\,\chi\br).
\ee
The first term in (\ref{a1}) can be written as
\[\g(\lambda) {}_1{\bf F}_2\bl(\!\!\begin{array}{c}\mu\!-\!k\\1\!-\!k, 1\!+\!\nu\end{array}\bl|\,\chi\br)
=\frac{\g(\lambda)}{\g(1+\nu)}\sum_{r=k}^\infty \frac{(\mu-k)_r \chi^r}{\g(1-k+r)(1+\nu)_r r!}\]
\[=\frac{\g(\lambda) \chi^k}{\g(1+\nu)}\sum_{n=0}^\infty\frac{(\mu-k)_{n+k}\chi^n}{(1+\nu)_{n+k} n! (n+k)!}
=\frac{\g(\mu) \chi^k}{k! \g(1+\nu+k)} \sum_{n=0}^\infty \frac{(\mu)_n \chi^n}{(k+1)_n(1+\nu+k)_n n!} \]
\[=\g(\mu)\chi^k {}_1{\bf F}_2\bl(\!\!\begin{array}{c}\mu\\k\!+\!1, 1\!+\!\nu\!+\!k\end{array}\bl|\,\chi\br).\]
Hence $Q(s)=0$ when $\mu-\lambda=k$, $k\geq 0$.

Next, consider the case $\mu-\lambda=-k$, with $k=1, 2, \ldots\ $. Then we have
\bee\label{a2}
Q(s)=\g(\lambda) {}_1{\bf F}_2\bl(\!\!\begin{array}{c}\mu\!+\!k\\k\!+\!1, 1\!+\!\nu\end{array}\bl|\,\chi\br)
- \g(\mu)\chi^{-k}{}_1{\bf F}_2\bl(\!\!\begin{array}{c}\mu\\1\!-\!k, 1\!+\!\nu\!-\!k\end{array}\bl|\,\chi\br).
\ee
The second term in (\ref{a2}) can be written in the form
\[\g(\mu)\chi^{-k}{}_1{\bf F}_2\bl(\!\!\begin{array}{c}\mu\\1\!-\!k, 1\!+\!\nu\!-\!k\end{array}\bl|\,\chi\br)
=\frac{\g(\mu)\chi^{-k}}{\g(1+\nu-k)}\sum_{r=k}^\infty \frac{(\mu)_r\chi^r}{\g(1-k+r)(1+\nu-k)_r r!}\]
\[=\frac{\g(\mu)}{\g(1+\nu-k)}\sum_{n=0}^\infty \frac{(\mu)_{n+k} \chi^n}{(1+\nu-k)_{n+k}n! (n+k)!}
=\frac{\g(\mu+k)}{k! \g(1+\nu)} \sum_{n=0}^\infty \frac{(\mu+k)_n \chi^n}{(1+\nu)_n (k+1)_n n!}\]
\[=\g(\lambda) {}_1{\bf F}_2\bl(\!\!\begin{array}{c}\mu\!+\!k\\k\!+\!1, 1\!+\!\nu\end{array}\bl|\,\chi\br)\br\}.\]
Hence $Q(s)=0$ when $\mu-\lambda=-k$, $k\geq 1$.

This concludes the demonstration that $H(s)$ is regular when $\mu-\lambda$ assumes integer values.

\vspace{0.6cm}

\begin{center}
{\bf Appendix B: Evaluation of the residue of the double pole at $s=1$}
\end{center}
\setcounter{section}{2}
\setcounter{equation}{0}
\renewcommand{\theequation}{\Alph{section}.\arabic{equation}}
When $\gamma+\nu=-1$ the function $H(s)$ in (\ref{e22}) has a double pole at $s=1$. We let $s=1+2\epsilon$, with $\epsilon\to 0$, so that $\lambda=\fs(s-1)=\epsilon$. Then, provided $\mu\neq 1, 2, \ldots\ $,
\[H(s)=\frac{\pi B}{\sin \pi(\mu-\epsilon)}\bl\{\g(\epsilon) {}_1{\bf F}_2\bl(\!\!\begin{array}{c}\epsilon\\1\!-\!\mu\!+\!\epsilon,1\!+\!\nu\end{array}\bl|\,\chi\br)-\g(\mu) \chi^{\mu-\epsilon}
{}_1{\bf F}_2\bl(\!\!\begin{array}{c}\mu\\1\!+\!\mu\!-\!\epsilon,1\!+\!\mu\!+\!\nu\!-\!\epsilon\end{array}\bl|\,\chi\br)\br\}\]
\[=\frac{\pi B}{\epsilon \sin \pi\mu}\frac{\{1+\epsilon \pi \cot \pi\mu+O(\epsilon^2)\}}{\g(1+\nu)}
\bl\{\frac{\g(1+\epsilon)}{\g(1-\mu+\epsilon)}{}_1 F_2\bl(\!\!\begin{array}{c}\epsilon\\1\!-\!\mu\!+\!\epsilon,1\!+\!\nu\end{array}\bl|\,\chi\br)\hspace{1cm}\]
\[\hspace{5cm}-\epsilon\g(\mu)\g(1+\nu) \chi^{\mu}{}_1{\bf F}_2\bl(\!\!\begin{array}{c}\mu\\1\!+\!\mu,1\!+\!\mu\!+\!\nu\end{array}\bl|\,\chi\br)+O(\epsilon^2)\br\}. \]
\[=\frac{B\g(\mu)}{\epsilon \g(1+\nu)} \{1+\epsilon \pi \cot \pi\mu+O(\epsilon^2)\}\bl\{(1-\epsilon({\hat\gamma}+\psi(1-\mu)) {}_1 F_2\bl(\!\!\begin{array}{c}\epsilon\\1\!-\!\mu\!+\!\epsilon,1\!+\!\nu\end{array}\bl|\,\chi\br)\]
\[\hspace{5cm}-\epsilon\g(\mu)\g(1+\nu) \chi^{\mu}{}_1{\bf F}_2\bl(\!\!\begin{array}{c}\mu\\1\!+\!\mu,1\!+\!\mu\!+\!\nu\end{array}\bl|\,\chi\br)+O(\epsilon^2)\br\}. \]

Now
\[\frac{\g(1+\epsilon)}{\g(1-\mu+\epsilon)}=\frac{1}{\g(1-\mu)}\{1-\epsilon({\hat\gamma}+\psi(1-\mu))+O(\epsilon^2)\}\]
and
\[{}_1 F_2\bl(\!\!\begin{array}{c}\epsilon\\1\!-\!\mu\!+\!\epsilon,1\!+\!\nu\end{array}\bl|\,\chi\br)=1+\frac{\epsilon\chi}{(1-\mu)(1+\nu)}\bl\{1+\frac{1! \chi}{(2-\mu)(2+\nu)2!}+\frac{2! \chi^2}{(2-\mu)_2(2+\nu)_2 3!}+\cdots\br\}+O(\epsilon^2)\]
\[=1+\epsilon\chi \g(1-\mu)\g(1+\nu)\,{}_2{\bf F}_3\bl(\!\!\begin{array}{c} 1,1\\2,2-\mu,1+\nu\end{array}\bl|\,\chi\br)+O(\epsilon^2).\]
Upon use of the identity $\psi(\mu)-\psi(1-\mu)=-\pi \cot \pi\mu$ we find
\[H(s)=\frac{B\g(\mu)}{\epsilon \g(1+\nu)}\{1+\epsilon (A-{\hat\gamma}-\psi(\mu))+O(\epsilon^2)\},\]
where
\bee\label{b0}
A=\g(1+\nu)\bl\{\chi \g(1-\mu)\,{}_2{\bf F}_3\bl(\!\!\begin{array}{c}1,1\\2,2\!-\!\mu, 2\!+\!\nu\end{array}\bl|\,\chi\br)-\frac{\pi\chi^\mu}{\sin \pi\mu}\,{}_1{\bf F}_2\bl(\!\!\begin{array}{c}\mu\\1\!+\!\mu, 1\!+\!\mu\!+\!\nu\end{array}\bl|\,\chi\br)\br\}.
\ee

Then, since $\zeta(1+2\epsilon)=(2\epsilon)^{-1}\{1+{\hat\gamma} 2\epsilon+O(\epsilon^2)\}$, we obtain the residue of $H(s)a^s \zeta(s)$ at $s=1$ when $\gamma+\nu=-1$ given by
\bee\label{b1}
\mbox{Res}_{s=1}=\frac{2aB\g(\mu)}{\g(1+\nu)}\bl\{\log\,a+\frac{1}{2}(A+{\hat\gamma}-\psi(\mu))\br\}\qquad (\mu\neq 1, 2,\ldots),
\ee
where we recall that $B=(b/2)^\nu/(2\g(\mu))$.

In the case of integer $\mu$, the limiting form of (\ref{b1}) is required. We consider here only the case $\mu=1$. Setting $\mu=1+\epsilon$, $\epsilon\to0$, we have from (\ref{b0})
\[A=-\frac{\chi}{\epsilon(1+\nu)} \,{}_2F_3\bl(\!\!\begin{array}{c}1,1\\2,1\!-\!\epsilon,2\!+\!\nu\end{array}\!\bl|\,\chi\br)+\frac{\chi^{1+\epsilon}\g(1+\nu)}{\epsilon\g(2+\epsilon)\g(2+\nu+\epsilon)}\,{}_1F_2\bl(\!\!\begin{array}{c}1\!+\!\epsilon\\2\!+\!\epsilon, 2\!+\!\nu\!+\!\epsilon\end{array}\!\bl|\,\chi\br)+O(\epsilon).\]
Since
\[(\alpha+\epsilon)_n=(\alpha)_n\{1+\epsilon\Delta\psi(\alpha+n)+O(\epsilon^2)\},\qquad\Delta\psi(\alpha+n):=\psi(\alpha+n)-\psi(\alpha),\]
we find that
\begin{eqnarray*}
{}_2 F_3\bl(\!\!\begin{array}{c}1,1\\2,1\!-\!\epsilon,2\!+\!\nu\end{array}\bl|\,\chi\br)&=&F+\epsilon\sum_{r=0}^\infty \frac{\chi^n}{(2)_n(2+\nu)_n}\, \Delta\psi(1+n)+O(\epsilon^2),\\
{}_1 F_2\bl(\!\!\begin{array}{c}1\!+\!\epsilon\\2\!+\!\epsilon,2\!+\!\nu\!+\!\epsilon\end{array}\bl|\,\chi\br)&=&
F+\epsilon\bl(\sum_{n=0}^\infty \frac{\chi^n}{(2)_n(2+\nu)_n}\, \Delta\psi(1+n)-F_*\br)+O(\epsilon^2),
\end{eqnarray*}
where $F:={}_1F_2(1;2,2+\nu;\chi)$ and
\bee\label{b3}
F_*:=\sum_{n=0}^\infty \frac{\chi^n}{(2)_n (2+\nu)_n} \,[\Delta\psi(2+n)+\Delta\psi(2+\nu+n)].
\ee

Then, after simplification,
\[A=-\frac{\chi F}{1+\nu} \{\psi(2)+\psi(2+\nu)-\log\,\chi\}-\frac{\chi F_*}{1+\nu}+O(\epsilon)\]
as $\epsilon\to0$.
Upon noting that $\chi F/(1+\nu)={\cal I}_\nu(b)-1$, where ${\cal I}_\nu(b)$ is defined in (\ref{e24b}), and that $\psi(1)=-{\hat\gamma}$, $\psi(2)=1-{\hat\gamma}$, we finally find
\bee\label{b2}
\mbox{Res}_{s=1}=\frac{aB}{\g(1+\nu)}\bl\{2\log\,a+2{\hat\gamma}+\kappa(1-{\cal I}_\nu(b))-\frac{\chi F_*}{1+\nu}\br\}
\ee
when $\mu=1$ and $\gamma+\nu=-1$, where $\kappa:=1-{\hat\gamma}+\psi(2+\nu)-\log\,\chi$.

\vspace{0.6cm}

\begin{center}
{\bf Appendix C: Inverse factorial expansion of ${}_1F_2(\chi)$ in (\ref{e31})}
\end{center}
\setcounter{section}{3}
\setcounter{equation}{0}
\renewcommand{\theequation}{\Alph{section}.\arabic{equation}}
We obtain the inverse factorial expansion of the hypergeometric function
\[F\equiv {}_1F_2\bl(\!\!\begin{array}{c}\!\!m-\fs s\\1\!-\!\mu\!+\!m\!-\!\fs s, 1\!+\!\nu\end{array}\bl|\,\chi\br), \qquad \chi:=\frac{b^2}{4}\]
appearing in (\ref{e31}) for large $s$ in $\Re (s)>0$. We write $F$ in the form
\[F=\sum_{n=0}^\infty \frac{a_n(s) \chi^n}{(1+\nu)_n n!},\]
where
\[a_n(s):=\frac{(m\!-\!\fs s)_n}{(1\!-\!\mu\!+\!m\!-\!\fs s)_n}=\frac{\g(1\!-\!m\!+\!\fs s)}{\g(1\!-\!m\!-\!n+\fs s)}\,\frac{\g(\mu\!-\!m\!-\!n\!+\!\fs s)}{\g(\mu\!-\!m\!+\!\fs s)}\]
\[=1+\frac{2(1-\mu)n}{s}+\frac{2(1-\mu)}{s^2}\{(\mu-2m)n+(\mu-2)n^2\}+O(s^{-3})\]
as $s\to\infty$ in $\Re (s)>0$,
upon use of the reflection formula for the gamma function and the asymptotic expansion of a ratio of two gamma functions given in \cite[p.~141]{DLMF}. This last result can then be arranged in the form
\bee\label{c1}
a_n(s)=1+\frac{A_1(n)}{s+\mu-1}+\frac{A_2(n)}{(s+\mu-1)(s+\mu-2)}+\cdots,
\ee where 
\[A_1(n)=2(1-\mu)n\qquad A_2(n)=2(1-\mu)\{(2m+1-\mu)n+(2-\mu)n(n-1)\}.\]

We therefore obtain 
\[F=\sum_{n=0}^\infty \frac{\chi^n}{(1+\nu)_n n!}\bl\{1+\frac{A_1(n)}{s\!+\!\mu\!-\!1}+\frac{A_2(n)}{(s\!+\!\mu\!-\!1)(s\!+\!\mu\!-\!2)}+\cdots \br\}\]
\[=\g(1+\nu)(\fs b)^{-\nu} I_\nu(b)+2(1-\mu)\bl\{\frac{1}{s\!+\!\mu\!-\!1}+\frac{2m\!+\!1\!-\!\mu}{(s\!+\!\mu\!-\!1)(s\!+\!\mu\!-\!2)}\br\} \sum_{n=1}^\infty\frac{n \chi^n}{(1+\nu)_n n!}\]
\[+\frac{2(1-\mu)(2-\mu)}{(s\!+\!\mu\!-\!1)(s\!+\!\mu\!-\!2)} \sum_{n=2}^\infty \frac{n(n-1) \chi^n}{(1+\nu)_n n!} +\cdots\ .\]
Now
\[\sum_{n=1}^\infty\frac{n \chi^n}{(1+\nu)_n n!}=\chi \sum_{n=0}^\infty \frac{\chi^n}{(1+\nu)_{n+1} n!}=\g(1+\nu) (\fs b)^{1-\nu} I_{\nu+1}(b)\]
and
\[\sum_{n=2}^\infty \frac{n(n-1) \chi^n}{(1+\nu)_n n!}=\chi^2 \sum_{n=0}^\infty \frac{\chi^n}{(1+\nu)_{n+2} n!}=\g(1+\nu) (\fs b)^{2-\nu} I_{\nu+2}(b),\]
where $I_\nu(b)$ denotes the modified Bessel function and we have used the fact that $(1+\nu)_{n+r}=(1+\nu+r)_n (1+\nu)_r$ for $r=1, 2, \ldots\ $.

Then the inverse factorial expansion for $F$ is given by
\bee\label{c2}
F=\g(1+\nu) (\fs b)^{-\nu} I_\nu(b) \bl\{C_0'+\frac{C_1'}{s\!+\!\mu\!-\!1}+\frac{C_2'}{(s\!+\!\mu\!-\!1)(s\!+\!\mu\!-\!2)}+\cdots\br\}
\ee
as $s\to\infty$ in $\Re (s)>0$, where
\[C_0'=1,\quad C_1'=(1-\mu)b\,\frac{I_{\nu+1}(b)}{I_\nu(b)},\]
\[ C_2'=(1-\mu)b\bl\{(2m\!+\!1\!-\!\mu)\,\frac{I_{\nu+1}(b)}{I_\nu(b)}+\frac{(2-\mu)b}{2}\,\frac{I_{\nu+2}(b)}{I_\nu(b)}\br\}.\]

\end{document}